\documentclass[12pt,psfig]{article}
\usepackage{graphicx,epsfig}
\usepackage{mathrsfs}
\usepackage{amssymb}
\def\bbb{\begin{eqnarray*}}

\def\eee{\end{eqnarray*}}

\begin{document}

\baselineskip=17pt

\begin{center}

\vspace{-0.6 in} {\large \bf  Invariance of deficiency indices of Hermitian subspaces
under relatively bounded perturbations$^*$}
\\ [0.3in]

Yan Liu $^{a}$, Yuming Shi $^{b, **}$\\

$^{a}$ Department of Mathematics and Physics, Hohai University,\\
Changzhou Campus 213022, P. R. China\\
$^{b}$ Department of Mathematics, Shandong University,\\
Jinan, Shandong 250100, P. R. China

\footnote{$^{*}$ This research was supported by the NNSF of China (Grant 11571202) and
the NSF of Jiangsu Province of China (Grant BK20170298).}
\footnote{$^{**}$ The corresponding author.}
\footnote{Email addresses: yanliumaths@126.com (Y. Liu), ymshi@sdu.edu.cn (Y. Shi)}
\end{center}

{\bf Abstract.}
This paper is concerned with the stability of deficiency indices of Hermitian subspaces (i.e., linear relations) under relatively bounded perturbations in Hilbert spaces.
Several results about invariance of deficiency indices of Hermitian subspaces under relatively bounded perturbations are established.
As a consequence, invariance of self-adjointness of Hermitian subspaces under relatively bounded perturbations is obtained.
In addition, it is shown that the deficiency indices may shrink in the special case that the relative bound is equal to 1.
The results obtained in the present paper generalize the corresponding results
for symmetric operators to more general Hermitian subspaces.
\medskip

\noindent{\bf 2010 AMS Classification}: 47A06, 47A55, 47B15.
\medskip

\medskip
\parindent=10pt

\section{ Introduction }

Perturbation theory occupies an important place in both pure and applied mathematics.
The perturbation theory of linear operators (i.e., single-valued operators) has been extensively studied
and a great deal of elegant results have been obtained (cf.,[2, 13, 16, 20, 24, 31, 32]).
In particular, stability of deficiency indices of symmetric operators under perturbations has received lots of attention.
We shall briefly recall several results about deficiency indices and their stability for symmetric operators.
In 1910, Weyl first studied the deficiency indices of second-order formally self-adjoint differential equations [33].
He showed that the deficiency index equals to the number of linearly independent square integrable solutions of the
differential equation for each spectral parameter $\lambda\in \mathbf{C}\backslash\mathbf{R}$.
Followed by this work, Atkinson studied the maximal deficiency index of Hamiltonian differential systems [2],
and obtained that the maximal deficiency index is invariant under bounded perturbation.
In particular, the stability of deficiency indices of symmetric operators under perturbations was deeply studied by Behncke and Focke [4].
They got the invariance of deficiency index of a symmetric operator under relatively bounded perturbation with relative bound less than 1.
This result extends the result for self-adjoint operators [16].
In the case that the relative bound of the perturbation is equal to 1, however, the deficiency index may shrink [4, Example].
For more results about the stability of the deficiency indices of symmetric operators,
we refer to [17, 36] and some references cited therein.

With further research of operator theory, more and more multi-valued operators and non-densely defined operators have been found.
For example, the operators generated by those linear continuous Hamiltonian systems,
which do not satisfy the definiteness conditions, and general linear symmetric difference systems
may be multi-valued or not densely defined in their corresponding Hilbert spaces (cf. [18, 23, 29]).
So the classical perturbation theory of linear operators is not available in this case.
Motivated by this need,
von Neumann [19] first introduced linear relations into functional analysis,
and then Arens [1] and many other scholars further studied and developed the fundamental theory of linear relations [7-12, 14, 15, 25, 27, 28].
A multi-valued operator is said to be a liner relation or a linear subspace (briefly, subspace)
since its graph is a linear subspace in its related product subspace.
We shall use the term ``subspace" in the present paper.

Since the fundamental theory of subspaces was established, the related perturbation problems have attracted extensive attention of scholars
and some good results have been obtained, including stability of boundness, closedness, self-adjointness, and spectrum of subspaces (cf. [3, 10, 21, 26, 30, 34, 35]).
However, compared to perturbation theory of linear operators, some important perturbation problems of subspaces have not been studied
or have been not thoroughly studied.
In the present paper, we are interested in stability of deficiency indices of Hermitian subspaces under perturbations.

It is well known that the deficiency indices of Hermitian subspaces play an important role in the theory of self-adjoint extensions.
According to the generalized von Neumann theory [7] and the GKN theory [25],
a Hermitian subspace has a self-adjoint extension if and only if its positive and negative deficiency indices are equal
and its self-adjoint extension domains have a close relationship with its deficiency indices.
So it is of great significance to determine the deficiency indices of Hermitian subspaces.

To the best of our knowledge, there seem a few results about the stability of the deficiency indices of Hermitian subspaces under perturbations,
though there are many results about the deficiency indices of Hermitian subspaces (cf.[2, 5, 6, 22, 23, 25]).
In 1998, Cross introduced concepts of boundedness, relative boundedness, and deficiency of linear relations,
and showed that the defciency of a subspace are stable under bounded perturbation under certain additional conditions [10, Corollary III.7.6].
The deficiency of a subspace has a close relationship with its deficiency indices.
In 2013, Zheng [37] obtained the invariance of the minimal and maximal deficiency indices
for discrete Hamiltonian systems under bounded perturbations.
In 2018, Ren [21] discussed the stability of index for linear relations,
and showed that the deficiency indices of dissipative linear relations are stable under relatively bounded perturbations with relative bounds less than 1 under certain additional conditions.
In the present paper, we shall focus on the study of the invariance of the deficiency indices of Hermitian subspaces under relatively bounded perturbations.
We shall remark that the results given in the present paper not only generalize the corresponding results given in [4]
for symmetric operators to Hermitian subspaces, but also cover the results obtained in [37] (see Remarks 3.3, 3.5, and 3.6).

The rest of the present paper is organized as follow.
In Section 2, some basic concepts and useful fundamental results about subspaces are introduced.
In particular, a useful lemma about subspaces is given,
which will take an important role in the study of stability of deficiency indices of Hermitian subspaces under perturbations.
In Section 3, stability of deficiency indices of Hermitian subspaces under perturbations is studied.
For any two Hermitian subspaces $S$ and $T$,
the deficiency indices of $T$ are invariant under the condition that $S$ is relatively bounded with respect to $T+tS$ for every $t\in[0,1]$.
This condition is weaker than that $S$ is relatively bounded with respect to $T$ with relative bound less than 1.
Therefore, as a consequence, the invariance of deficiency indices of Hermitian subspaces is obtained
under relatively bounded perturbations with relative bounds less than 1.
In addition, it is shown that the deficiency indices may shrink in the case that the relative bound is equal to 1.
As a consequence, stability of self-adjointness of Hermitian subspaces under relatively bounded perturbations is obtained (see Corollary 3.2).
\medskip

\noindent\textit{\bf Remark 1.1.}
We shall apply the results obtained in the present paper to
discuss the invariance of deficiency indices of second-order symmetric linear difference equations and
discrete linear Hamiltonian systems under relatively bounded perturbations, separately, in our forthcoming papers.
\medskip

\section{Preliminaries}

In this section, we shall recall some basic concepts, introduce some fundamental results,
and give some results about subspaces, which will be used in the sequent section.

By ${\mathbf C}$ and ${\mathbf R}$ denote the sets of the complex numbers and the real numbers, respectively.
Let $X$ be a complex Hilbert space with inner product $\langle\cdot,\cdot\rangle$,
and $T$ a linear subspace (briefly, subspace) in the product space $X^2:=X\times X$.
The domain $D(T)$, range $R(T)$, and null space $N(T)$ of $T$ are respectively defined by
\vspace{-0.2cm}
$$\begin{array}{rrll}
D(T):&=&\{x\in X:\, (x,f)\in T \;{\rm for\; some}\;f\in X\},\\[0.4ex]
R(T):&=&\{f\in X:\, (x,f)\in T \;{\rm for\; some}\;x\in X\},\\[0.4ex]
N(T):&=&\{x\in X:\, (x,0)\in T \}.
\end{array}\vspace{-0.2cm}$$
Further, denote
$$T(x):=\{f\in X:\, (x,f)\in T \},\;\;T^{-1}:=\{(f,x):\, (x,f)\in T \}.$$
It is evident that $T(0)=\{0\}$ if and only if $T$ can uniquely determine a single-valued linear
operator from $D(T)$ into $X$ whose graph is $T$.
For convenience, a linear operator in $X$ will always be identified with a subspace in $X^2$ via its graph.

Let $T$ and $S$ be two subspaces in $X^2$ and $\alpha \in {\mathbf C}$.
Define\vspace{-0.2cm}
$$\begin{array}{cl}
\alpha T:=\{(x,\alpha f):\, (x,f)\in T\},\\[0.6ex]
T+S:= \{(x,f+g):\,(x,f)\in T, (x,g)\in S\},\\[0.6ex]
ST:=\{(x,g):\, (x,f)\in T, (f,g)\in S \:{\rm for\: some}\: f\in X\}.\\[0.6ex]
\end{array}\vspace{-0.2cm}$$

A subspace $T$ is called closed if it is a closed subspace in $X^2$.
It is evident that if $T$ is closed, then $T-\alpha I_{id}$ is closed,
where $I_{id}:=\{(x, x):\, x\in X\}$. Without any confusion, we briefly denote it by $I$.

The adjoint of $T$ is defined by\vspace{-0.2cm}
$$T^*=\{(y,g)\in X^2:\,\langle g,x\rangle=\langle y,f\rangle\;{\rm
for\;all}\; (x,f)\in T\}.\vspace{-0.2cm}$$
A subspace $T\subset X^2$ is called a Hermitian subspace if $T\subset T^*$,
and it is called a self-adjoint subspace if $T=T^*$.
In addition, $T$ is a Hermitian subspace if and only if
$\langle f, y\rangle=\langle x, g\rangle$ for all $(x,f), (y,g)\in T$.
\medskip

\noindent{\bf Definition 2.1 {\rm [25, Definition 2.3]}.}  Let $T$ be a subspace in $X^2$ and $\lambda \in \mathbf{C}$.
The subspace $R(T-\lambda I)^\perp$ is called the deficiency space of $T$ and $\lambda$,
and the number $d_\lambda(T):=\dim(R(T-\lambda I))^\perp$ is called
the deficiency index of $T$ and $\lambda$.
\medskip

It can be easily verified that the deficiency indices of $T$ and its closure with the same $\lambda$ are equal.
Let $T$ be a Hermitian subspace in $X^2$. By [25, Theorem 2.3], $d_\lambda(T)$ is constant
in the upper and lower half-planes; that is, $d_\lambda(T)=d_+(T)$ for all $\lambda\in \mathbf{C}$ with ${\rm Im}\lambda>0$
and $d_\lambda(T)=d_-(T)$ for all $\lambda\in \mathbf{C}$ with ${\rm Im}\lambda<0$, where $d_{\pm}(T):=d_{\pm i}(T)$.
The pair $(d_+(T), d_-(T))$ is called the deficiency indices of $T$, and $d_+(T)$ and $d_-(T)$ are
called the positive and negative deficiency indices of $T$, respectively.
\medskip

\noindent\textit{{\bf Lemma 2.1 {\rm [10, Proposition I.2.8 and Corollary I.2.11]}.}  Let $T$ be a subspace in $X^2$.
\begin{itemize}\vspace{-0.2cm}
\item[{\rm (i)}] For any given $x\in D(T)$, $y\in T(x)$ if and only if $T(x)=\{y\}+T(0)$.
In particular, $0\in T(x)$ if and only if $T(x)=T(0)$;\vspace{-0.2cm}
\item[{\rm (ii)}] $(TT^{-1})(y)=\{y\}+T(0)$ for every $y\in R(T)$ and $(T^{-1}T)(x)=\{x\}+T^{-1}(0)$ for every $x\in D(T)$.\vspace{-0.2cm}
\end{itemize}}
\medskip

In the following, we shall recall concepts of the norm of a subspace and relatively boundedness of two subspaces, and their fundamental properties.

Let $E$ be a closed subspace of $X$. Define the following quotient space [16]:
$$X/E:=\{[x]:\,x\in X\},\quad [x]=\{x\}+E.$$
We define an inner product on the quotient space $X/E$ by
$$\langle [x],[y] \rangle:=\langle x^\bot, y^\bot\rangle,\quad  [x],[y]\in X/E,                                \eqno(2.1)$$
where $x=x_0+x^\bot$, $y=y_0+y^\bot$ with $x_0,\,y_0\in E$ and $x^\bot,\,y^\bot\in E^\bot.$
It can been easily verified that the above inner product is well-defined and $X/E$ with this inner product is a Hilbert space.
The norm induced by this inner product is the same as that of $X/E$ induced by the norm of $X$.

Now, define the following natural quotient map:
$$Q_E^X:\:X \rightarrow X/E,\:x \mapsto [x].$$
Let $T$ be a subspace in $X^2$.
By $Q_T$ or simply $Q$ denote $Q_{\overline{T(0)}}^X$ for briefness without confusion. Define
$$\tilde{T}_s=G(Q_T)T.                                                                                     \eqno(2.2)$$
Then $\tilde{T}_s$ is a linear operator with domain $D(T)$ [10, Proposition II.1.2].
The norm of $T$ at $x\in D(T)$ and the norm of $T$ are defined by, respectively (see [10, II.1]),
$$\begin{array}{rrll}
&&\|T(x)\|:=\|\tilde{T}_s(x)\|,\\[0.6ex]
&&\|T\|:=\|\tilde{T}_s\|=\sup\{\|\tilde{T}_s(x)\|:\,x\in D(T)\, with\, \|x\|\leq 1\}.\\[0.6ex]
\end{array}                                                                                                 \eqno(2.3)$$

\noindent\textit{{\bf Lemma 2.2 {\rm [10, Propositions II.1.4-II.1.7]}.}
Let $T$ and $S$ be two subspaces in $X^2$. Then
\begin{itemize}\vspace{-0.2cm}
\item[{\rm (i)}] $\|T(x)\|=d(T(x),0)=d(T(x),T(0))=d(y,\overline{T(0)})=d(y,T(0))$ for every $x\in D(T)$ and $y\in T(x)$;\vspace{-0.2cm}
\item[{\rm (ii)}] $\|(\alpha T)(x)\|=|\alpha|\|T(x)\|$,\quad $\|\alpha T\|=|\alpha|\|T\|$ for every $x\in D(T)$ and $\alpha\in \mathbf{C}$;\vspace{-0.2cm}
\item[{\rm (iii)}] $\|S(x)+T(x)\| \leq \|S(x)\|+\|T(x)\|$ for $x\in D(T)\cap D(S)$, \quad $\|S+T\| \leq \|S\|+\|T\|$.\vspace{-0.2cm}
\end{itemize}}
\medskip

\noindent{\bf Definition 2.2 {\rm [10, Definition VII.2.1]}.} Let $T$ and $S$ be two subspaces in $X^2$.
\begin{itemize}\vspace{-0.2cm}
\item [{\rm (1)}] The subspace $S$ is said to be $T$-bounded if $D(T)\subset D(S)$ and
there exists a constant $c\geq 0$ such that
$$\|S(x)\|\leq c(\|x\|+\|T(x)\|),\quad x\in D(T).                            \vspace{-0.2cm}$$
\item [{\rm (2)}] If $S$ is $T$-bounded, then the infimum of all numbers $b\geq 0$
for which a constant $a\geq 0$ exists such that
$$\|S(x)\|\leq a \|x\|+b\|T(x)\|,\quad x\in D(T),                                                                \eqno(2.4)$$
is called the $T$-bound of $S$.\vspace{-0.2cm}
\end{itemize}

\noindent\textit{\bf Remark 2.1.}
Condition (2.4) is equivalent to the following condition:
$$\|S(x)\|^2\leq a'^2\|x\|^2+b'^2\|T(x)\|^2,\quad x\in D(T),                                                    \eqno(2.5)$$
where the constants $a'\geq 0$ and $b'\geq 0$. In fact, it can be easily deduced that (2.5) implies (2.4) with $a=a', b=b'$, whereas (2.4)
implies (2.5) with $a'^2=(1+\varepsilon^{-1})a^2$ and $b'^2=(1+\varepsilon)b^2$ with an arbitrary $\varepsilon>0$.
Consequently, the $T$-bound of $S$ may as well be defined as the infimum of the possible values of $b'$.
This equivalency is also be pointed out in [35, Lemma 5.10] without proof.
Here, we give its detailed statement.
\medskip

\noindent\textit{{\bf Lemma 2.3 {\rm [30, Theorem 6.3]}.}
Let $S$ and $T$ be subspaces in $X^2$ with $D(T)\subset D(S)$ and $S(0)\subset T(0)$.
Assume that $S$ is $T$-bounded with $T$-bound less than $1$.
Then $T+S$ is closed if and only if $T$ is closed.}\medskip

\noindent\textit{{\bf Lemma 2.4 {\rm [30, Proposition 2.1]}.}
Let $S$ and $T$ be two subspaces in $X^2$.
Then $T=T-S+S$ if and only if $D(T)\subset D(S)$ and $S(0)\subset T(0)$.}\medskip

\noindent\textit{{\bf Lemma 2.5 {\rm [35, Lemma 5.8]}.} Let $T$ be a self-adjoint subspace in $X^2$.
If $S$ is a Hermitian subspace in $X^2$ and $D(T)\subset D(S)$, then  $S(0)\subset T(0)$.}\medskip

\noindent\textit{{\bf Lemma 2.6 {\rm [30, Propositions 3.1 and 3.3]}.}
Let $T$ be a Hermitian subspaces in $X^2$.
Then $D(T)\subset {T(0)}^\perp$ and
$$\langle \tilde{T}_s(x_2),[x_1] \rangle=\langle [x_2], \tilde{T}_s(x_1) \rangle,\,\,x_1, x_2\in D(T).$$}
\vspace{-0.5cm}

Now, we present two well-known results about operators, which are useful in the sequent discussions.\medskip

\noindent\textit{{\bf Lemma 2.7 {\rm [31, Theorem 4.23]}.} If $P$ and $Q$ are orthogonal projections on $X$
such that $\|P-Q\|<1$, then $\dim R(P)=\dim R(Q)$ and $\dim R(I-P)=\dim R(I-Q)$.}
\medskip

\noindent\textit{{\bf Lemma 2.8 {\rm [25, Lemma 2.14]}.}
Let $T:\,X\rightarrow X$ be a linear operator, $Y$ be a subspace in $X$, and $R(T)\subset \overline{Y}$.
Then, for each $x\in D(T)$,
$$\|T(x)\|=\sup\{|\langle T(x), y\rangle|:\,y\in Y\, {\rm with} \,\|y\|=1\}.$$}
\vspace{-0.5cm}

The following result will take an important role in the study
of stability of deficiency indices of Hermitian subspaces under perturbations.
It is a generalization of Theorem 5.25 in [31] for operators to subspaces.\medskip

\noindent\textit{{\bf Lemma 2.9.} Let $A$ and $B$ be two subspaces in $X^2$ with $D(A)\subset D(B)$ and $B(0)\subset A(0)$.
Assume that there exists a constant $c\geq 0$ such that
$$\|B(x)\|\leq c\|A(x)\|, \:x\in D(A).                                                                              \eqno(2.6)$$
For ever $k\in \mathbf{C}$, let $P_k$ denote the orthogonal projection from $X$ onto $\overline{R(A+kB)}$.
Then $\|P_k-P_0\|\rightarrow 0$ as $k\rightarrow 0$.}\medskip

\noindent{\bf Proof.}
In the case that $c=0$, the result holds obviously.

Now, we consider the case that $c>0$.
We shall show that the result holds by [31, Theorem 4.33].
Since $D(A)\subset D(B)$ and $B(0)\subset A(0)$, we get that $A=A+kB-kB$ by Lemma 2.4 for every $k\in \mathbb{C}$ with $|k|\leq {1/(2c)}$.
Then, for all $x\in D(A)$, it follows from (2.6) that
$$\|B(x)\|\leq c\|A(x)\|\leq c\|(A+kB)(x)\|+c|k|\|B(x)\|\leq c\|(A+kB)(x)\|+\frac{1}{2}\|B(x)\|,$$
which implies that
$$\|B(x)\|\leq 2c\|(A+kB)(x)\|,\: x\in D(A).                                                                        \eqno(2.7)$$

Fix any $h\in R(P_0)^\perp=R(A)^\perp$.
Since $X=\overline{R(A+kB)}\oplus R(A+kB)^\perp$,
there exist $h_1\in\overline{R(A+kB)}$ and $h_2\in R(A+kB)^\perp$ such that $h=h_1+h_2$ and $P_k(h)=h_1$.
For any $(y,g)\in A+kB$, there exist $(y,g_1)\in A$ and $(y,g_2)\in B$ such that $g=g_1+kg_2$.
Then
$$\langle P_k(h),g\rangle=\langle h_1,g\rangle=\langle h,g\rangle=\langle h,g_1+kg_2\rangle=\bar{k}\langle h,g_2\rangle.\eqno(2.8)$$
By $B(0)\subset A(0)$ we have that $A(0)^\bot \subset B(0)^\bot$.
Thus, $h\in R(A)^\perp \subset A(0)^\perp \subset B(0)^\perp$.
In addition, there exist $g_{2,0}\in \overline{B(0)}$ and $g_2^\perp \in B(0)^\perp$ such that $g_2=g_{2,0}+g_2^\perp$.
So it follows from (2.1) that
$$\langle h,g_2\rangle=\langle h,g_{2,0}+g_2^\perp \rangle=\langle h,g_2^\perp \rangle=\langle [h],[g_{2}] \rangle=\langle [h],\tilde{B}_s(y) \rangle,$$
where $[h],\,[g_{2}]\in X/\overline{B(0)}$. This, together with (2.7) and (2.8), yields that
$$\begin{array}{rrll}
&&|\langle P_k(h),g\rangle| \leq |k|\|[h]\|\|\tilde{B}_s(y)\|=|k|\|h\|\|B(y)\|\\[0.5ex]
&\leq & 2c|k|\|h\|\|(A+kB)(y)\| \leq 2c|k|\|h\|\|g\|,
\end{array}\vspace{-0.2cm}                                                                                    \eqno(2.9)$$
in which (2.1) and (i) of Lemma 2.2 have been used.
By Lemma 2.8, it follows from (2.9) that
$$\|P_k(h)\|\leq 2c|k|\|h\|.                                                                                     \eqno(2.10)$$

On the other hand, fix any $h'\in R(P_k)^\perp=R(A+kB)^\perp$.
There exist $h'_1\in\overline{R(A)}$ and $h'_2\in R(A)^\perp$ such that $h'=h'_1+h'_2$ and $P_0(h')=h'_1$.
Then, for any $(x,f)\in A$ we have that
$$\langle P_0(h'),f\rangle=\langle h'_1,f\rangle=\langle h',f\rangle.                                         \eqno(2.11)$$
By the assumption that $D(A)\subset D(B)$, there exists $g$ such that $(x,g)\in B$.
It is evident that $f+kg\in R(A+kB)$. So, by (2.11) we have that
$$\langle P_0(h'),f\rangle=\langle h',f+kg-kg\rangle=\langle h',f+kg\rangle-\bar{k}\langle h',g\rangle=-\bar{k}\langle h',g\rangle. \eqno(2.12)$$
Since $(A+kB)(0)=A(0)\subset R(A+kB)$, $h'\in R(A+kB)^\perp \subset A(0)^\perp \subset B(0)^\perp$.
In addition, there exist $g_{0}\in \overline{B(0)}$ and $g^\perp \in B(0)^\perp$ such that $g=g_{0}+g^\perp$.
Again by (2.1) we get that
$$\langle h',g\rangle=\langle h',g_{0}+g^\perp \rangle=\langle h',g^\perp \rangle=\langle [h'],[g] \rangle=\langle [h'],\tilde{B}_s(x) \rangle,$$
where $[h'],\,[g]\in X/\overline{B(0)}$.
This, together with (2.12) and (2.6), implies that
$$|\langle P_0 (h'),f\rangle| \leq |k|\|[h']\|\|\tilde{B}_s(x)\|
=|k|\|h'\|\|B(x)\| \leq c|k|\|h'\|\|A(x)\| \leq c|k|\|h'\|\|f\|,                                             \eqno(2.13)$$
in which (2.1) and (i) of Lemma 2.2 have been used.
Again by Lemma 2.8, it follows from (2.13) that
$$\|P_0(h')\|\leq c|k|\|h'\|.                                                                                  \eqno(2.14)$$
Therefore, by [31, Theory 4.33], (2.10), and (2.14) we get that $\|P_k-P_0\|\leq 2c|k|$,
which implies that $\|P_k-P_0\|\rightarrow 0$ as $k\rightarrow 0$.
The proof is complete.
\medskip

\noindent\textit{\bf Remark 2.2.}
Lemma 2.9 is a generalization of [31, Theorem 5.25] for operators to subspaces.
In 2012, we gave its another generalization (see Lemma 2.15 in [25]).
We shall remark that Lemma 2.9 weakens the conditions of Lemma 2.15 in [25],
in which it is required that $\|g\|\leq c\|f\|$ for all $(x,f)\in A$ and $(x,g)\in B$ with some $c\geq 0$.
\medskip

The following result comes from [10, (a) of Theorem III 4.2] or [30, Theorem 5.1].\medskip

\noindent\textit{{\bf Lemma 2.10.}
Let $T$ be a closed subspace in $X^2$. Then $T$ is bounded if and only if $D(T)$ is closed.}
\medskip

\section{Invariance of deficiency indices of Hermitian subspaces under perturbations}

In this section, we shall study the stability of deficiency indices of Hermitian subspaces under perturbations.
\medskip

We shall first prove the following two useful lemmas.\medskip

\noindent\textit{{\bf Lemma 3.1.} Let $T$ be a Hermitian subspace in $X^2$. Then
\begin{itemize}\vspace{-0.2cm}
\item [{\rm (i)}] for any $z=a+ib\in \mathbf{C}$ with $a,\,b\in\mathbf{R}$,
$$\|(T-zI)(x)\|^2=\|(T-aI)(x)\|^2+b^2\|x\|^2,\,x\in D(T);$$
\item [{\rm (ii)}] for any $z\in \mathbf{C} \backslash \mathbf{R}$,
$\|(T-zI)^{-1}\|\leq |{\rm Im} z|^{-1}$.\vspace{-0.2cm}
\end{itemize}}
\medskip

\noindent{\bf Proof.}
(i) Fix any $x\in D(T)$ and any $z=a+ib\in \mathbf{C}$ with $a,\,b\in\mathbf{R}$. We have that
$$\begin{array}{rrll}
&&\|(T-zI)(x)\|^2=\|Q_T(T-zI)(x)\|^2\\[1ex]
&=&\|Q_TT(x)-aQ_T(x)-ibQ_T(x)\|^2\\[1ex]
&=&\langle \tilde{T}_s(x)-a[x]-ib[x],\tilde{T}_s(x)-a[x]-ib[x]\rangle,\\[1ex]
\end{array}                                                                                                \eqno(3.1)$$
where $[x]\in X/{\overline{T(0)}}$.
Since $T$ is Hermitian, it follows from Lemma 2.6 that $x\in T(0)^\perp$ and
$$\langle \tilde{T}_s(x),[x] \rangle=\langle [x],\tilde{T}_s(x) \rangle.$$
Hence, we obtain that
$$\langle \tilde{T}_s(x)-a[x],-ib[x] \rangle=\langle ib[x],\tilde{T}_s(x)-a[x] \rangle.$$
Inserting the above relation into (3.1), we get that
$$\|(T-zI)(x)\|^2=\|\tilde{T}_s(x)-a[x]\|^2+b^2\|[x]\|^2=\|(T-aI)(x)\|^2+b^2\|[x]\|^2.                     \eqno(3.2)$$
In addition, $\|[x]\|=\|x\|$ by noting that $x\in T(0)^\perp$.
Therefore, it follows from (3.2) that Assertion (i) holds.

(ii) Fix any $z=a+ib$ with $a,\,b\in\mathbf{R}$ and $b\neq 0$.
For any $(x,y)\in T-zI$, by Assertion (i) one has that
$$\|(T-zI)(x)\|^2=\|(T-aI)(x)\|^2+b^2\|x\|^2\geq b^2\|x\|^2,$$
which yields that
$$\|x\|\leq \frac{1}{|b|}\|(T-zI)(x)\|.$$
This, together with Lemma 2.2, yields that
$$\|(T-zI)^{-1}(y)\|\leq \|x\|\leq \frac{1}{|b|}\|(T-zI)(x)\|\leq \frac{1}{|b|}\|y\|,                       \eqno(3.3) $$
which yields that
$$\|(T-zI)^{-1}\|\leq \frac{1}{|b|}=\frac{1}{|{\rm Im} z|}.$$
The whole proof is complete.
\medskip

\noindent\textit{\bf Remark 3.1.}
Assertion (i) of Lemma 3.1 extends the result given in [16, p. 270] for closed symmetric operators to Hermitian subspaces.
Assertion (ii) of Lemma 3.1 extends [11, Theorem 3.20] for self-adjoint subspaces and [27, Theorem 3.4] for closed Hermitian subspaces to Hermitian subspaces.
\medskip

\noindent\textit{{\bf Lemma 3.2.}
Let $T$ and $S$ be two subspaces in $X^2$ with $D(T)\subset D(S)$ and $S(0)\subset T(0)$.
If $S$ is $T$-bounded with $T$-bound less than 1,
then $S$ is $(T+tS)$-bounded for every $t\in [0,1]$.}
\medskip

\noindent{\bf Proof.}
Fix any $t\in [0,1]$.
By the assumption that $S$ is $T$-bounded with $T$-bound less than 1,
there exist $a\geq 0$ and $0\leq b<1$ such that
$$\|S(x)\|\leq a\|x\|+b\|T(x)\|, \,\,x\in D(T).                                                               \eqno(3.4)$$
Since $D(T)\subset D(S)$ and $S(0)\subset T(0)$, we have that $T=T+tS-tS$ by Lemma 2.4.
Thus, for any $x\in D(T)$, it follows from (3.4) that
$$\begin{array}{rrll}
\|S(x)\|
&\leq& a\|x\|+b\|T(x)\|\\[0.6ex]
&\leq& a\|x\|+b\|(T+tS)(x)\|+bt\|S(x)\|\\[0.6ex]
&\leq& a\|x\|+b\|(T+tS)(x)\|+b\|S(x)\|,
\end{array}$$
which implies that
$$\|S(x)\|\leq \frac{a}{1-b}\|x\|+\frac{b}{1-b}\|(T+tS)(x)\|.$$
Therefore, $S$ is $(T+tS)$-bound.
This completes the proof.
\medskip

\noindent\textit{\bf Remark 3.2.} We shall point out that ${b}/{(1-b)}$ is not necessarily less than 1 in Lemma 3.2.
So the $(T+tS)$-bound of $S$ is not necessarily less than 1 for every $t\in [0,1]$.
\medskip

The following theorem is the main result of the present paper, in which the condition is weaker than
relative boundedness with relative bound less than 1.
\medskip

\noindent\textit{{\bf Theorem 3.1.} Let $T$ and $S$ be Hermitian subspaces in $X^2$ with $D(T)\subset D(S)$ and $S(0)\subset T(0)$.
If $S$ is $(T+tS)$-bounded for every $t\in [0,1]$, then $d_\pm(T+S)= d_\pm(T)$.}
\medskip

\noindent{\bf Proof.}
Fix any $t\in [0,1]$.
By the assumption that $S$ is $(T+tS)$-bounded, there exist $a>0$ and $b>0$ such that
$$\|S(x)\|^2\leq a^2\|x\|^2+b^2\|(T+tS)(x)\|^2, \,\,x\in D(T).                                                          \eqno(3.5)$$
Let $z=\pm i{a}/{b}$. Note that $T+tS$ is Hermitian.
It follows from Lemma 3.1 that
$$\|(T+tS-zI)(x)\|^2=\|(T+tS)(x)\|^2+\frac{a^2}{b^2}\|x\|^2, \,\,x\in D(T),$$
which, together with (3.5), yields that
$$\|S(x)\|^2\leq b^2\|(T+tS-zI)(x)\|^2, \,\,x\in D(T).$$
Hence, we get that
$$\|S(x)\|\leq b\|(T+tS-zI)(x)\|, \,\,x\in D(T).$$
For every $k\in \mathbb{C}$, let $P_{t,k}$ denote the orthogonal projection from $X$ onto $\overline{R(T+tS-zI+kS)}$.
By Lemma 2.9 we get that
$$\|P_{t,k}-P_{t,0}\|\rightarrow 0 \:{\rm as}\: k\rightarrow 0.                                                         \eqno(3.6)$$
Let $F(t)=I-P_{t,0}$, $t\in [0,1]$. Then the operator-valued function $F(t)$ is continuous in $[0,1]$ by (3.6).
By Lemma 2.7 and the Henine-Borel theorem, we get that ${\rm{dim}}R(F(0))={\rm{dim}}R(F(1))$;
that is, ${\rm{dim}}R(T-zI)^\perp={\rm{dim}}R(T+S-zI)^\perp$, and so $d_\pm(T+S)= d_\pm(T)$.
The proof is complete.
\medskip

The following result is a direct consequence of Lemma 3.2 and Theorem 3.1.
\medskip

\noindent\textit{{\bf Corollary 3.1.}
Let $T$ and $S$ be Hermitian subspaces in $X^2$ with $D(T)\subset D(S)$ and $S(0)\subset T(0)$.
If $S$ is $T$-bounded with $T$-bound less than 1, then $d_\pm(T+S)= d_\pm(T)$.}
\medskip

\noindent\textit{\bf Remark 3.3.}
Corollary 3.1 extends [32, Theorem 6.2] or the first Corollary 2 in [4] for closed symmetric operators to Hermitian subspaces.
\medskip

Combining Corollary 3.1 and Lemma 2.3, we can get the following result:\medskip

\noindent\textit{{\bf Corollary 3.2.}
Let $T$ and $S$ be Hermitian subspaces in $X^2$ with $D(T)\subset D(S)$ and $S(0)\subset T(0)$.
If $S$ is $T$-bounded with $T$-bound less than 1, then $T+S$ is self-adjoint if and only if $T$ is self-adjoint.}
\medskip

\noindent\textit{\bf Remark 3.4.} The sufficiency of Corollary 3.2 still holds when the assumption that $S(0)\subset T(0)$ is dropped.
In fact, if $T$ is self-adjoint, then the condition that $D(T)\subset D(S)$ implies that $S(0)\subset T(0)$ by Lemma 2.5.
The sufficiency of Corollary 3.2 is the same as that of [35, Theorem 5.2].
Here, Corollary 3.2 shows that the condition is not only sufficient but also necessary.
\medskip

\noindent\textit{{\bf Corollary 3.3.}
Let $T$ and $S$ be Hermitian subspaces in $X^2$ with $D(T)\subset D(S)$ and $S(0)\subset T(0)$.
If $N(T)\subset N(S)$ and $ST^{-1}$ is bounded with bound less than 1, then $d_\pm(T+S)= d_\pm(T)$.}
\medskip
\vspace{-0.5cm}

\noindent{\bf Proof.}
For any $x\in D(T)$, by Lemma 2.1 we have that
$$(T^{-1}T)(x)=\{x\}+T^{-1}(0)=\{x\}+N(T).$$
Hence, $(ST^{-1}T)(x)=S(x)$ by $D(T)\subset D(S)$ and $N(T)\subset N(S)$.
Consequently, by [10, (3) of Remarks II.3.14] one has that
$$\|S(x)\|=\|(ST^{-1}T)(x)\|\leq \|ST^{-1}\| \|T(x)\|.$$
Note that $ST^{-1}$ is bounded with bound less than 1.
It follows that $S$ is $T$-bounded with $T$-bound less than 1.
Therefore, $d_\pm(T+S)= d_\pm(T)$ by Corollary 3.1.
This completes the proof.
\medskip

\noindent\textit{{\bf Corollary 3.4.}
Let $T$ be a Hermitian subspace in $X^2$ and $S$ a Hermitian operator in $X$ with $D(T)\subset D(S)$.
If $N(T)\subset N(S)$, then $ST^{-1}$ is an operator in $X$.
In addition, if $ST^{-1}$ is bounded with bound less than 1, then $d_\pm(T+S)= d_\pm(T)$.}
\medskip

\noindent{\bf Proof.}
By Corollary 3.3, it suffices to show that $ST^{-1}$ is an operator in $X$.
For any $(0,g)\in ST^{-1}$, there exists $x\in X$ such that $(x,0)\in T$, $(x,g)\in S$.
Thus, $x\in N(T)\subset N(S)$, and so $g=S(x)=0$ by the assumption that $S$ is an operator.
This implies that $ST^{-1}$ is an operator in $X$.
The proof is complete.
\medskip

\noindent\textit{{\bf Corollary 3.5.}
Let $T$ and $S$ be Hermitian subspaces in $X^2$ with $D(T)\subset D(S)$ and $S(0)\subset T(0)$.
If $S$ is $T$-bounded and satisfies that
$${\rm{Re}}\langle f,g \rangle \geq 0 ,\,\,\, (x,f)\in T,\,(x,g)\in S,                                            \eqno(3.7)$$
then $d_\pm(T+S)= d_\pm(T)$.}
\medskip

\noindent{\bf Proof.}
By Theorem 3.1, it suffices to show that $S$ is $(T+tS)$-bounded for every $t\in [0,1]$.
By the assumption that $S$ is $T$-bounded and Remark 2.1, there exist nonnegative numbers $a$ and $b$ such that
$$\|S(x)\|^2\leq a^2\|x\|^2+b^2\|T(x)\|^2, \,\,x\in D(T).                                                             \eqno(3.8)$$
For any $(x,f)\in T$ and $(x,g)\in S$,
it follows from (3.7) that for any $t\in [0,1]$,
$$\|f+tg\|^2= \|f\|^2+t^2\|g\|^2+2t{\rm{Re}}\langle f,g \rangle \geq \|f\|^2+t^2\|g\|^2\geq \|T(x)\|^2+t^2\|S(x)\|^2,$$
which, together with (3.8), yields that
$$b^2\|f+tg\|^2+a^2\|x\|^2 \geq b^2\|T(x)\|^2+b^2t^2\|S(x)\|^2+a^2\|x\|^2\geq (1+b^2t^2)\|S(x)\|^2.$$
Hence,
$$\|S(x)\|\leq \frac{b}{\sqrt{1+b^2t^2}}\|f+tg\|+\frac{a}{\sqrt{1+b^2t^2}}\|x\|,\,x\in D(T).$$
By the arbitrariness of $(x,f)\in T$ and $(x,g)\in S$, it follows that
$$\|S(x)\|\leq \frac{b}{\sqrt{1+b^2t^2}}\|(T+tS)(x)\|+\frac{a}{\sqrt{1+b^2t^2}}\|x\|,\,x\in D(T).$$
This implies that $S$ is $(T+tS)$-bounded for every $t\in [0,1]$.
This completes the proof.
\medskip

\noindent\textit{{\bf Corollary 3.6.}
Let $T$ and $V$ be two Hermitian subspaces in $X^2$ satisfying $D(T)=D(V)=D$ and $V(0)=T(0)$.
If there exist constants $a\geq 0$ and $0\leq b<1$ such that
$$\|(V-T)(x)\|\leq a\|x\|+b(\|T(x)\|+\|V(x)\|),\,\,x\in D,                                                     \eqno(3.9)$$
then $d_\pm(V)= d_\pm(T)$. }
\medskip

\noindent{\bf Proof.}
Let $S=V-T$ and fix any $t\in [0,1]$.
Then $D(S)=D$ and $S(0)=T(0)=V(0)$.
It follows from Lemma 2.4 that $T=T+tS-tS$ for every $t\in [0,1]$ and $V=T+S$.
This, together with (3.9), yields that
$$\begin{array}{rrll}
\|S(x)\|&\leq& a\|x\|+b {\Big(} \|T(x)+tS(x)-tS(x)\|+\|T(x)+tS(x)+(1-t)S(x)\| {\Big )}\\[0.6ex]
&\leq& a\|x\|+2b\|(T+tS)(x)\|+b\|S(x)\|,\,x\in D,\,t\in [0,1],\\[0.6ex]
\end{array}$$
which implies that
$$\|S(x)\|\leq \frac{a}{1-b}\|x\|+\frac{2b}{1-b}\|(T+tS)(x)\|,\,x\in D,\,t\in [0,1].$$
Hence, $S$ is $T+tS$-bounded for every $t\in [0,1]$. Therefore, by Theorem 3.1 we get that $d_\pm(T+S)= d_\pm(T)$.
And consequently, $d_\pm(V)= d_\pm(T)$.
The proof is complete.
\medskip

In Corollary 3.1, the assumption that the relative bound of $S$ with respect to $T$ is less than 1 cannot be dropped in general.
The following result is concerned with the case that the relative bound is equal to 1.
It will be shown that the deficiency index may shrink in this case in general (see Remark 3.5).\medskip

\noindent\textit{{\bf Theorem 3.2.}
Let $T$ be a closed Hermitian subspace in $X^2$ and $S$ a symmetric operator (densely defined and Hermitian operator) in $X$ with $D(T)\subset D(S)$.
If there exists $a\geq 0$ such that
$$\|S(x)\|\leq a\|x\|+\|T(x)\|, \,\,x\in D(T),                                                                       \eqno(3.10)$$
then $d_\pm(T+S)\leq d_\pm(T)$.}
\medskip

\noindent{\bf Proof.}
We first show that $d_+(T+S)\leq d_+(T)$.
Without loss of generality we assume that $d_+(T)=:m<\infty$.
Let $X_t:=R(T+tS-iI)$ for $t\in [0,1]$.

Fix any $t\in [0,1)$.
Then $tS$ is $T$-bounded with $T$-bound less than 1 by (3.10).
Since $T$ is closed, it follows from Lemma 2.3 that $T+tS$ is closed.
Thus, $T+tS-iI$ and $(T+tS-iI)^{-1}$ are closed.
Since $T+tS$ is Hermitian, it follows from (ii) of Lemma 3.1 that $(T+tS-iI)^{-1}$ is bounded.
Hence, it follows from Lemma 2.10 that $X_t$ is closed.
Therefore, $X=X_t\oplus X_t^\perp$ for every $t\in [0,1)$.
In addition, by Corollary 3.1 we have that $d_+(T)=d_+(T+tS)={\rm{dim}}X_t^\perp=m<\infty$ for every $t\in [0,1)$.

Suppose that $d_+(T+S)={\rm{dim}}X_1^\perp=:n<+\infty$ and $m<n$.
Then there exists a normalized element $y_t\in X_t\bigcap X_1^\perp$ for every $t\in [0,1)$.
Choose a subsequence $\{y_{t_n}\}_{n=1}^{\infty}$ with $t_n\in [0,1)$ for each $n \geq 1$
satisfying that $t_n\rightarrow 1$ and $y_{t_n}\rightarrow y$ as $n\rightarrow \infty$.
Then $y\in X_1^\perp$ and $\|y\|=1$.
There exists $x_{t_n}\in D(T)$ such that $(x_{t_n},y_{t_n})\in T+t_nS-iI$ for each $n \geq 1$.
With a similar argument to that used for (3.3), one has that $\|x_{t_n}\|\leq \|y_{t_n}\|=1$ for all $n \geq 1$.
Thus, $\{x_{t_n}\}_{n=1}^{\infty}$ is bounded.
Let $h_{t_n}:=S(x_{t_n})$ for each $n \geq 1$. Then, $(x_{t_n},y_{t_n}-t_nh_{t_n}+ix_{t_n})\in T$ for each $n \geq 1$.
By (3.10) and Lemma 2.2 we get that
$$\|h_{t_n}\| \leq a\|x_{t_n}\|+\|T(x_{t_n})\| \leq a\|x_{t_n}\|+\|y_{t_n}-t_nh_{t_n}+ix_{t_n}\| \leq a+2+t_n\|h_{t_n}\|,\,n \geq 1,$$
which yields that
$$\|(1-t_n)h_{t_n}\|\leq a+2,\,n \geq 1.$$
Thus, $\{(1-t_n)h_{t_n}\}_{n=1}^{\infty}$ is bounded.
Note that every bounded set in $X$ is weakly sequentially compact by [13, Theorems I.6.15 and I.7.21].
So there exist a subsequence $\{(1-t_{n_{k}})h_{t_{n_{k}}}\}_{k=1}^{\infty}$ and an $w\in X$
such that $\{(1-t_{n_{k}})h_{t_{n_{k}}}\}_{k=1}^{\infty}$ converges weakly to $w$.

For any $u\in D(S)$, since $S$ is a symmetric operator, we get that
$$\langle u,w\rangle=\lim_{k\to\infty}\langle u,(1-t_{n_{k}})h_{t_{n_{k}}}\rangle=
\lim_{k\to\infty}\langle u,(1-t_{n_{k}})S(x_{t_{n_{k}}})\rangle=\lim_{k\to\infty}(1-t_{n_{k}})\langle S(u),x_{t_{n_{k}}}\rangle=0.$$
This, together with the fact that $S$ is densely defined, yields that $\langle u,w\rangle=0$ holds for every $u\in X$.
It follows that
$$\lim_{k\to\infty}\langle y,y_{t_{n_{k}}}+(1-t_{n_{k}})h_{t_{n_{k}}}\rangle=\|y\|^2+\langle y,w\rangle=1.                                                \eqno(3.11)$$
By noting that $(x_{t_{n_{k}}},y_{t_{n_{k}}}+(1-t_{n_{k}})h_{t_{n_{k}}})\in T+S-iI$,
it follows that $y_{t_{n_{k}}}+(1-t_{n_{k}})h_{t_{n_{k}}}\in X_1$.
Thus, $\langle y,y_{t_{n_{k}}}+(1-t_{n_{k}})h_{t_{n_{k}}}\rangle=0$ by the fact that $y\in X_1^\perp$.
This contradicts (3.11).
Hence, $d_+(T+S)\leq d_+(T)$.

In the case that $d_+(T+S)=+\infty$, choose an $(m+1)$-dimensional subspace $M\subset X_1^\perp$.
With a similar argument to that used in the case that $d_+(T+S)<+\infty$, one can get a contradiction by replacing $X_1^\perp$ with $M$ in this case.

With a similar argument to the above, one can prove that $d_{-}(T+S)\leq d_{-}(T)$.
The entire proof is complete.
\medskip

As a consequence, the following result can be directly derived from Theorem 3.2.\medskip

\noindent\textit{{\bf Corollary 3.7.}
Let $T$ be a self-adjoint subspace in $X^2$ and $S$ a symmetric operator in $X$ with $D(T)\subset D(S)$.
If there exists $a\geq 0$ such that (3.10) holds,
then $T+S$ is essential self-adjoint.}
\medskip

\noindent\textit{\bf Remark 3.5.}
\begin{itemize}\vspace{-0.2cm}
\item [{\rm (1)}]
We shall point out that in the case that the $T$-bound of $S$ is equal to 1, the result of Theorem 3.2 is optimal.
For example, let $T$ be a self-adjoint operator and $S=-T$.
Then (3.10) holds and $d_\pm(T+S)= d_\pm(T)=0$.
Consider another example:
Let $T$ be a closed symmetric operator with $d_+(T)=d_-(T)>0$ and $S=(T^*T)^{{1}/{2}}$.
Then, $D(T)=D(S)$, $\|T(x)\|=\|S(x)\|$ for all $x\in D(T)$, and $S$ is a self-adjoint operator by [31, p. 198].
Hence, (3.10) holds and $d_\pm(T+S)=0$ by Corollary 3.7. Consequently, $d_\pm(T+S)<d_\pm(T)$.
\item [{\rm (2)}] Theorem 3.2 is a generalization of the second Corollary 1 in [4] for closed symmetric operators to Hermitian subspaces.
\end{itemize}

\medskip

\noindent\textit{{\bf Corollary 3.8.}
Let $T$ and $S$ be Hermitian subspace in $X^2$ with $D(T)\subset D(S)$ and $S(0)\subset T(0)$.
If there exists $a\geq 0$ such that (3.10) holds and $S$ is $(T+S)$-bounded,
then $d_\pm(T+S)= d_\pm(T)$.}
\medskip

\noindent{\bf Proof.}
It follows from (3.10) that for any $t\in [0,1)$ and $x\in D(T)$,
$$\|S(x)\| \leq a\|x\|+\|(T+tS)(x)-tS(x)\|
\leq  a\|x\|+\|(T+tS)(x)\|+t\|S(x)\|
$$
which yields that
$$\|S(x)\| \leq \frac{a}{1-t}\|x\|+\frac{1}{1-t}\|(T+tS)(x)\|.$$
This implies that $S$ is $(T+tS)$-bound for every $t\in [0,1)$.
This, together with the assumption that $S$ is $(T+S)$-bounded, yields that
the conditions in Theorem 3.1 hold, and so $d_\pm(T+S)= d_\pm(T)$ by Theorem 3.1.
The proof is complete.
\medskip

\noindent\textit{\bf Remark 3.6.}
By the definition of relative boundedness for subspaces,
we shall remark that the results about stability of deficiency indices obtained in the present paper still hold under bounded perturbations.
Therefore, the results in this paper cover the results given in [37, Theorems 3.1 and 3.2], in which the stability of the minimal and maximal deficiency indices
for discrete Hamiltonian systems under bounded perturbations was studied.
\medskip

\bigskip

%
%

\end{document}